%
%
%
%
%
%
\magnification=\magstep1
\input amstex
\documentstyle{amsppt}
\voffset-2pc
\nologo

\define\C{{\Bbb C}}

\define\dee{\partial}
\redefine\O{\Omega}
\redefine\phi{\varphi}
\define\Obar{\overline{\Omega}}
\define\Ot{\widetilde\Omega}
\define\Oh{\widehat\Omega}

\NoRunningHeads
\topmatter
\title
Finitely generated function fields and complexity in potential
theory in the plane
\endtitle
\author Steven R. Bell${}^*$ \endauthor
\thanks ${}^*$Research supported by NSF grant DMS-9623098 \endthanks
\address
Mathematics Department, Purdue University, West Lafayette, IN  47907 USA
\endaddress
\email bell\@math.purdue.edu \endemail
\abstract
We prove that the Bergman kernel function
associated to a finitely connected domain $\O$ in the plane is given
as a rational combination of only {\it three\/} basic functions of one complex
variable: an Alhfors map, its derivative, and one other function whose
existence is deduced by means of the field of meromorphic functions on
the double of $\O$.  Because many other functions of conformal mapping
and potential theory can be expressed in terms of the Bergman kernel,
our results shed light on the complexity of these objects.  We also
prove that the Bergman kernel is an {\it algebraic\/} function of a single
Ahlfors map and its derivative.  It follows that many objects of potential
theory associated to a multiply connected domain are algebraic if and only
if the domain is a finite branched cover of the unit disc via an algebraic
holomorphic mapping.
\endabstract
\endtopmatter

\document

\hyphenation{bi-hol-o-mor-phic}
\hyphenation{hol-o-mor-phic}

\subhead 1. Introduction \endsubhead
On a simply connected domain in the plane, the Riemann mapping
function can be used to pull back the classical kernel functions
of potential theory from the unit disc and it follows that the
kernel functions are given as simple rational combinations of
two functions of one complex variable, the Riemann map and its
derivative.  Indeed, if $a$ is a point in a simply connected domain
$\O\ne\C$ and $f_a(z)$ is the Riemann mapping function mapping
$\O$ one-to-one onto the unit disc with $f_a(a)=0$ and $f_a'(a)>0$,
then the Bergman kernel $K(z,w)$ associated to $\O$ is given by
$$K(z,w)=
\frac{f_a'(z)\overline{f_a'(w)}}
{\pi (1-f_a(z)\overline{f_a(w)})^2}.$$
Another way to write the same formula is
$$K(z,w)=
\frac{c K(z,a)\overline{K(w,a)}}
{(1-f_a(z)\overline{f_a(w)})^2}$$
where $c=\pi/f_a'(a)^2$.  In this paper, I shall prove that,
similarly, on an $n$-connected domain $\O$ such that no boundary
component is a point, the Bergman kernel is a rational combination
of only {\it three\/} basic functions of {\it one\/} complex variable.
One of the basic functions is an Ahlfors mapping associated to the
domain.  We shall define Ahlfors maps carefully in the next section.
Suffice it say for now that the Ahlfors mapping $f_a$ associated to a
point $a\in\O$ is a branched $n$-to-one covering map of $\O$ onto the
unit disc with $f_a(a)=0$ and with $f_a'(a)>0$ and maximal.  In many
ways, the Ahlfors map can be thought of as a``Riemann mapping function''
for a multiply connected domain.  The following theorem shows that
the Ahlfors map takes over the role of the Riemann map in the kernel
identities mentioned above.

\proclaim{Theorem 1.1}
Suppose that $\O$ is an $n$-connected domain in the plane such that no
boundary component is a point.  There exist points $a$ and $b$ in $\O$
with the property that the Bergman kernel $K(z,w)$ associated
to $\O$ is a rational combination of the Ahlfors map $f_a$, its
derivative $f_a'$, and the function of one variable given by $K(\cdot,b)$.
To be precise, there exists a formula of the form
$$K(z,w)=
R(f_a(z),f_a'(z),K(z,b),
\,\overline{f_a(w)}\,,
\,\overline{f_a'(w)}\,,
\,\overline{K(w,b)}\,)\tag1.1$$
where $R$ is a {\it rational\/} function on $\C^6$.
\endproclaim

There is nothing particularly special about the points $a$ and $b$ in
Theorem~1.1.  In fact, the proof shall show that there is a dense open
subset of points $(a,b)$ in $\O\times\O$ satisfying the same properties.
We shall prove Theorem~1.1 in \S\S2-4.

I shall also prove that, under the hypotheses of Theorem~1.1, there
is an irreducible polynomial of two complex variables $P(z,w)$ such
that the two functions $f_a(z)$ and $K(z,b)/f_a'(z)$ satisfy 
$P(f_a(z),K(z,b)/f_a'(z))\equiv0$, and so it follows that
$K(z,b)$ is an algebraic function of $f_a$ and $f_a'$.  Consequently,
formula (1.1) yields that $K(z,w)$ is an algebraic function
of $f_a(z)$, $f_a'(z)$, and  the conjugates of $f_a(w)$ and $f_a'(w)$.
The polynomial $P(z,w)$ is a rather interesting algebraic geometric
object attached to $\O$.  We show that it is conformally invariant,
but we do not explore its possibly deeper significance here.

The Bergman kernel is at the heart of the Bergman metric and Ahlfors maps
determine the Carath\'eodory metric.  The results of this paper show
that the two metrics are connected in an algebraic, but perhaps very
complicated, manner on multiply connected domains.

Because most other objects of potential theory can be expressed in
terms of the Bergman kernel, the results of this paper shall
yield that many of the objects of potential theory on a finitely
connected domain are given as rational combinations
of an Ahlfors map, its derivative, and {\it one\/}
other basic function of one variable.  For example the classical
``functions of the first kind,'' $F_1'$, $F_2'$,\dots,$F_{n-1}'$,
associated to an $n$-connected domain will be seen to be rational
combinations of the three basic functions mentioned above, and it will
follow that the square of the Szeg\H o kernel is a rational combination
of the three basic functions.  Furthermore, the relationships proved
in \cite{4} between the Szeg\H o kernel, the Poisson kernel, and the
Green's function will shed light on the complexity of all these classical
objects of potential theory.

An interesting example that is easily analyzed by the methods of this
paper is the $2$-connected domain $A(r)$ given by $\{z\,:\, |z+1/z|<r\}$
where $r$ is a real constant bigger than $2$.  We shall show that
the Bergman and Szeg\H o kernels associated to $A(r)$ are algebraic and
that every $2$-connected domain in the plane such that no boundary
component is a point is biholomorphic to exactly one domain $A(r)$ with $r>2$.
We shall show that biholomorphic maps between domains with algebraic
Bergman kernel functions must be algebraic, and it shall
follow that we may think of $A(r)$ as being the defining member of a conformal
class.  We shall also show that there is a polynomial $P(z,w)$ of two
complex variables such that the Bergman kernel associated to $A(r)$
satisfies
$$P(K(z,w),f(z),\,\overline{f(w)}\,)=0$$
on $A(r)\times A(r)$ where $f(z)=(1/r)(z+1/z)$.
Let $\hat A(r)$ denote the double of $A(r)$.  The polynomial equation
will reveal that it is possible to analytically continue $K(z,w)$
to $\hat A(r)\times\hat A(r)$ as a finitely valued multivalued function
with algebraic singularities.

Before we proceed to give all the details, I shall sketch the proof of the
Theorem~1.1.  For the purposes of this introduction,
suppose that $\O$ is a bounded $n$-connected domain in the plane
bounded by $n$ non-intersecting real analytic curves.  The first step
in the proof is to show that the Bergman kernel generates itself in
the following sense.  Let $K(z,w)$ denote the Bergman kernel associated
to $\O$.  Let $K_0(z,w)$ also denote the Bergman kernel and let $K_m(z,w)$
denote the function $(\dee^m/\dee\bar w^m)K(z,w)$.  We prove that there
exists a finite subset $\Cal A$ of $\O$ and a positive integer
$N$ such that the Bergman kernel $K(z,w)$ is given as a rational
combination of functions from the finite set of functions of $z$,
$$\{ K_m(z,b) \,:\,  b\in\Cal A,\ 0\le m\le N\},$$
and the finite set of functions of $w$,
$$\{ \overline{K_m(w,b)} \,:\,  b\in\Cal A,\ 0\le m\le N\}.$$
Let $f_a(z)$ denote an Ahlfors mapping associated to $\O$ which
maps $a\in\O$ to the origin.  We next show that $f_a(z)$ and functions
of $z$ of the form $$\frac{K_m(z,b)}{f_a'(z)}$$
extend to the double of $\O$ as meromorphic functions.
Next, we show that there are points $a$ and $b$ in $\O$ such
that $f_a(z)$ and $K(z,b)/f_a'(z)$ form a primitive pair for
the field of meromorphic functions on the double of $\O$.  This means
that any meromorphic function on the double of $\O$ can be written
as a rational combination of these two functions (see
Farkas and Kra \cite{10, page~249}).  Finally, it follows that the
functions of $z$ of the form
$$\frac{K_m(z,w)}{f_a'(z)}$$
can be expressed as rational combinations of
$f_a(z)$ and $K(z,b)/f_a'(z)$.  Since there are finitely many
functions of $z$ of the form $K_m(z,w)$ that generate the
Bergman kernel, the existence of formula (1.1) will be proved.

Theorem 1.1 can be interpreted to mean that the Ahlfors maps
play the same role in the multiply connected setting that the
Riemann maps play for simply connected domains.  It is surprising,
however, that the same reasoning as above will show that {\it any\/}
proper holomorphic mapping of the domain onto the unit disc can be
used in place of an Ahlfors map --- even one of very high order.

I proved in \cite{4} that the Bergman and Szeg\H o kernels associated
to an $n$-connected domain are generated by $n+1$ basic functions of
one complex variable.  It is ironic that in order to prove that the
Bergman kernel is generated by only the {\it three\/} basic functions of one
variable mentioned above, I shall need to prove as an intermediate step
that it is generated by $3(n^2-2n+2)$ functions of one variable.

\subhead 2. Background information and the statement of a main theorem
\endsubhead
Before we can begin to carefully state and prove our main results, we
need to review some known facts about the classical kernel functions.
Many of these facts and formulas can be found in Stefan Bergman's book
\cite{8}.  I have also written up most of these results in \cite{2} in
the same spirit as this paper and I include cross references here to give
the interested reader access to a uniform approach to the whole subject.

To begin with, we shall assume that $\O$ is a bounded $n$-connected
domain in the plane with $C^\infty$ smooth boundary.  (Later, we shall
consider general $n$-connected domains such that no boundary component
is a point.)

Let $\gamma_j$, $j=1,\dots,n$,
denote the $n$ non-intersecting $C^\infty$ simple closed curves which
define the boundary $b\O$ of $\O$, and suppose that $\gamma_j$ is
parameterized in the standard sense by $z_j(t)$, $0\le t\le 1$.  We
shall use the convention that $\gamma_n$ denotes the {\it outer
boundary curve\/} of $\O$.  Let $T(z)$ be the $C^\infty$ function
defined on $b\O$ such that $T(z)$ is the complex number representing
the unit tangent vector at $z\in b\O$ pointing in the direction of
the standard orientation (meaning that $iT(z)$ represents the {\it inward
pointing normal vector\/} at $z\in b\O$).  This complex unit tangent vector
function is characterized by the equation $T(z_j(t))=z_j'(t)/|z_j'(t)|$.

The symbol $A^\infty(\O)$ will denote the space of holomorphic functions
on $\O$ that are in $C^\infty(\Obar)$.  The space of complex valued
functions on $\O$ that are square integrable with respect to Lebesgue
area measure $dA$  will be denoted by $L^2(\O)$, and the space of complex
valued functions on $b\O$ that are square integrable with respect to arc
length measure $ds$ by $L^2(b\O)$.  The Bergman space of
holomorphic functions on $\O$ that are in $L^2(\O)$ will be denoted by
$H^2(\O)$ and the Hardy space of functions in $L^2(b\O)$ that are the
$L^2$ boundary values of holomorphic functions on $\O$ by
$H^2(b\O)$.  The Bergman projection $P_B$ is the orthogonal projection
of $L^2(\O)$ onto $H^2(\O)$, and the Szeg\H o projection $P_S$ is the
orthogonal projection of $L^2(b\O)$ onto $H^2(b\O)$.  The Bergman kernel
$K(z,w)$ and the Szeg\H o kernel $S(z,w)$ are the kernels for the respective
projections in the sense that
$$\gather
(P_B\phi)(z)=
\iint_{w\in\O}K(z,w)\,\phi(w)\  dA, \\
(P_S\psi)(z)=
\int_{w\in b\O}S(z,w)\,\psi(w)\  ds.
\endgather$$

The Bergman kernel $K(z,w)$ is related to the Szeg\H o kernel via the
identity
$$K(z,w)=4\pi S(z,w)^2+\sum_{i,j=1}^{n-1}
A_{ij}F_i'(z)\overline{F_j'(w)},\tag 2.1$$
where the functions $F_i'(z)$ are classical functions of potential theory
described as follows (\cite{8, page~119}, or see also
\cite{2, pages~94--96}).  The harmonic function $\omega_j$ which solves the
Dirichlet problem on $\O$ with boundary data equal to one on the boundary
curve $\gamma_j$ and zero on $\gamma_k$ if $k\ne j$ has a multivalued
harmonic conjugate.  The function $F_j'(z)$ is a globally defined
single valued holomorphic function on $\O$ which is locally defined as
the derivative of $\omega_j+iv$ where $v$ is a local harmonic
conjugate for $\omega_j$.  The Cauchy-Riemann equations reveal that
$F_j'(z)=2(\dee\omega_j/\dee z)$.  A very important fact that we shall
need is that the matrix $[A_{ij}]$ appearing in formula (2.1) is
non-singular.  That this is so was proved by Hejhal in \cite{11} (see
Lemma~1 on page 74 and Theorems 30 and 31 on pages~81--83).

The Bergman and Szeg\H o kernels are holomorphic in the first variable
and antiholomorphic in the second on $\O\times\O$ and they are hermitian,
i.e.,  $K(w,z)=\overline{K(z,w)}$.  Furthermore, the Bergman and
Szeg\H o kernels are in
$C^\infty((\Obar\times\Obar)-\{(z,z):z\in b\O\})$ as functions of $(z,w)$
(see \cite{2, page~100}).

We shall also need to mention the
Garabedian kernel $L(z,w)$, which is related to the Szeg\H o
kernel via the identity
$$\frac{1}{i} L(z,a)T(z)=S(a,z)\qquad\text{for $z\in b\O$ and $a\in\O$.}
\tag2.2$$
For fixed $a\in\O$, the kernel $L(z,a)$ is a holomorphic function of $z$
on $\O-\{a\}$ with a simple pole at $a$ with residue $1/(2\pi)$.
Furthermore, as a function of $z$, $L(z,a)$ extends to the boundary
and is in the space $C^\infty(\Obar-\{a\})$.  In fact, $L(z,w)$
is in $C^\infty((\Obar\times\Obar)-\{(z,z):z\in\Obar\})$ as a function
of $(z,w)$ (see \cite{2, page~102}).  Also, $L(z,a)$ is non-zero for all
$(z,a)$ in $\Obar\times\O$ with $z\ne a$ and $L(a,z)=-L(z,a)$ (see
\cite{2, page~49}).

For each point $a\in\O$, the function of $z$ given by
$S(z,a)$ has exactly $(n-1)$ zeroes in $\O$ (counting multiplicities) and
does not vanish at any points $z$ in the boundary of $\O$ (see
\cite{2, page~49}).

Given a point $a\in\O$, the Ahlfors map $f_a$ associated to the pair $(\O,a)$
is a proper holomorphic mapping of $\O$ onto the unit disc.  It is an
$n$-to-one mapping (counting multiplicities), it extends to be in
$A^\infty(\O)$, and it maps each boundary curve $\gamma_j$ one-to-one
onto the unit circle.  Furthermore, $f_a(a)=0$, and $f_a$ is the unique
function mapping $\O$ into the unit disc maximizing the quantity $|f_a'(a)|$
with $f_a'(a)>0$.  The Ahlfors map is related to the Szeg\H o kernel
and Garabedian kernel via (see \cite{2, page~49})
$$f_a(z)=\frac{S(z,a)}{L(z,a)}.\tag2.3$$
Note that $f_a'(a)=2\pi S(a,a)\ne 0$.  Because $f_a$ is $n$-to-one, $f_a$
has $n$ zeroes.  The simple pole of $L(z,a)$ at $a$ accounts for the simple
zero of $f_a$ at $a$.   The other $n-1$ zeroes of $f_a$ are given by the
$(n-1)$ zeroes of $S(z,a)$ in $\O-\{a\}$.  Let $a_1,a_2,\dots,a_{n-1}$ denote
these $n-1$ zeroes (counted with multiplicity) and let
$\Cal Z(a)$ denote the set $\{a,a_1,\dots,a_{n-1}\}$
which is equal to the set of zeroes of $f_a$ on $\O$.  Next, we define
a set $\Cal A(a)= \Cal Z(a)\cup_{k=1}^{n-1}\Cal Z(a_k)$.  I proved in \cite{3}
(see also \cite{2, page~105}) that, as $a$ tends to a boundary curve
$\gamma_j$, the $n-1$ zeroes $a_1,\dots,a_{n-1}$ become distinct simple zeroes
which separate and tend toward the $n-1$ distinct boundary components of $\O$
which are different from $\gamma_j$.  To be precise, as $z$ approaches
a point in $\gamma_j$, there is exactly one zero of $S(z,a)$ that approaches
a point in $\gamma_k$, $k\ne j$.  It
follows from this result that there is a finite subset $G$ of $\O$
such that $S(z,a)$ has $n-1$ distinct simple zeroes in $\O$ as a function
of $z$ for every point $a$ in $\O-G$.  Hence, for points $a$ in $\O-G$
that are sufficiently close to the boundary, the set ${\Cal A}(a)$ is the
set consisting of $a$ together with the simple zeroes $a_1,\dots,a_{n-1}$
together with the $n-1$ simple zeroes of $S(z,a_k)$ for each $k$,
$1\le k\le n-1$.  Note that the set ${\Cal A}(a)$ has at
most $n+(n-1)(n-2)=n^2-2n+2$ elements (because $a$ is one of the
zeroes of $S(z,a_k)$ for each $k$).  One of the key ingredients in
the work that follows is that the Bergman kernel finitely generates
itself.  To make precise what we mean by this, we must define some
function spaces.  Let $K_0(z,w)$ and $K(z,w)$ both denote the
Bergman kernel associated to $\O$ and let $K_m(z,w)$ denote the
function $(\dee^m/\dee\bar w^m)K(z,w)$.
Let $\Cal R_z({\Cal A}(a),N)$ denote the field of holomorphic functions
of $z$ that are finite rational combinations of
$$\{K_m(z,\alpha)\,:\, \alpha\in\Cal A(a)\text{ and }0\le m\le N\}$$
where $N$ is a positive integer.  Let
$\Cal R_{\bar w}({\Cal A}(a),N)$ denote the field of antiholomorphic
functions of $w$ that are finite rational combinations of
$$\{\overline{K_m(w,\alpha)}\,:\,\alpha\in\Cal A(a)\text{ and }0\le m\le N\},$$
and let
$\Cal R_{z,\bar w}({\Cal A}(a),N)$ denote the field of functions of $z$ and
$w$ that are finite rational combinations of functions in
$\Cal R_z({\Cal A}(a), N)$ and
$\Cal R_{\bar w}({\Cal A}(a),N)$.

\proclaim{Theorem 2.1}
Suppose that $\O$ is an $n$-connected domain in the plane such that no
boundary component of $\O$ is a point.  For points $a\in\O$ that are
sufficiently close to the boundary, the Bergman kernel function $K(z,w)$
associated to $\O$ is in $\Cal R_{z,\bar w}({\Cal A}(a),2)$.
Furthermore, $S(z,w)^2$, is in $\Cal R_{z,\bar w}({\Cal A}(a),2)$ and
the functions $F_j'(z)$, $j=1,\dots,n-1$ are each in $\Cal R_z({\Cal A}(a),2)$.
\endproclaim

We remark that similar reasoning can be used to prove that
the kernels $\Lambda(z,w)$ and $L(z,w)^2$ are in
$\Cal R_{z,w}({\Cal A}(a),2)$ (where the missing bar over the
$w$ means that the generating functions in $w$ are taken to
be holomorphic instead of conjugate-holomorphic).

The proof of Theorem 2.1 will be given in \S3.

Fix a point $a$ in $\O$ so that the zeroes $a_1,\dots,a_{n-1}$ of
$S(z,a)$ are distinct simple zeroes.  I proved in \cite{4, Theorem~3.1}
that the Szeg\H o kernel can be expressed in terms of the $n+1$
functions of one variable, $S(z,a)$, $f_a(z)$, and $S(z,a_i)$, $i=1,\dots,n-1$
via the formula
$$S(z,w)=\frac{1}{1-f_a(z)\overline{f_a(w)}}\left(c_0
S(z,a)\overline{S(w,a)}+
\sum_{i,j=1}^{n-1} c_{ij}S(z,a_i)\,\overline{S(w,a_j)}\right)\tag2.4$$
where $f_a(z)$ denotes the Ahlfors map associated to $(\O,a)$, $c_0=1/S(a,a)$,
and the coefficients $c_{ij}$ are given as the
coefficients of the inverse matrix to the matrix $\left[S(a_j,a_k)\right]$.

A similar identity exists for the Garabedian kernel (see \cite{4}).
$$L(z,w)=\frac{f_a(w)}{f_a(z)-f_a(w)}\left(c_0 S(z,a)L(w,a)+
\sum_{i,j=1}^{n-1} c_{ij}S(z,a_i)L(w,a_j)\right)\tag2.5$$
where the constants $c_0$ and $c_{ij}$ are the same as the constants
in (2.4).

Let $\Cal F'$ denote the vector space of functions given by the complex
linear span of the set of functions $\{F_j'(z):j=1,\dots,n-1\}$
mentioned above.  It is
a classical fact that $\Cal F'$ is $n-1$ dimensional.  It shall be
important for us to relate the functions in $\Cal F'$ to the Szeg\H o
and Bergman kernel functions.
Notice that $S(z,a_i)L(z,a)$ is in $A^\infty(\O)$ because the pole
of $L(z,a)$ at $z=a$ is cancelled by the zero of $S(z,a_i)$ at $z=a$.
Similarly, $S(z,a)L(z,a_i)$ is in $A^\infty(\O)$ because the pole
of $L(z,a_i)$ at $z=a_i$ is cancelled by the zero of $S(z,a)$ at $z=a_i$.
A theorem due to Schiffer (\cite{14}, see also \cite{2, page~80}) states
that the set of $n-1$ functions
$\{S(z,a_i)L(z,a):i=1,\dots,n-1\}$ form a basis for $\Cal F'$.
It is also shown in \cite{2, page~80} that the linear span of
$\{S(z,a_i)L(z,a):i=1,\dots,n-1\}$ is the same as the linear span
of $\{L(z,a_i)S(z,a):i=1,\dots,n-1\}$.  Hence, formula (2.1) can
be rewritten in the form
$$K(z,w)=4\pi S(z,w)^2+\sum_{i,j=1}^{n-1}
\lambda_{ij}{\Cal L}_i(z)\,\overline{{\Cal L}_j(w)}\tag2.6$$
where
${\Cal L}_i(z)=L(z,a_i)S(z,a)$ and, because the matrix
$[A_{ij}]$ in formula (2.1) is non-singular, the change of basis
we have used yields a matrix $[\lambda_{ij}]$ that is also non-singular.

The Bergman kernel is related to the classical Green's function via
(\cite{8, page~62}, see also \cite{2, page~131})
$$K(z,w)=-\frac{2}{\pi}\frac{\dee^2 G(z,w)}{\dee z\dee\bar w}.$$
Another kernel function on $\O\times\O$ that we shall need is given by
$$\Lambda(z,w)=-\frac{2}{\pi}\frac{\dee^2 G(z,w)}{\dee z\dee w}.$$
In the literature, this function is sometimes written as $L(z,w)$ with
anywhere between zero and three tildes and/or hats over the top.
We have chosen the symbol $\Lambda$ here to avoid
confusion with our notation for the Garabedian kernel above.
It follows from well known properties of the Green's function that
$\Lambda(z,w)$ is holomorphic in $z$ and $w$  and is in
$C^\infty(\Obar\times\Obar-\{(z,z):z\in\Obar\})$.
If $a\in\O$, then $\Lambda(z,a)$ has a double pole at
$z=a$ as a function of $z$ and $\Lambda(z,a)=\Lambda(a,z)$ (see
\cite{2, page~134}).

The Bergman kernel is related to $\Lambda$ via the identity
$$\Lambda(w,z)T(z)=-K(w,z)\overline{T(z)}\qquad\text{for $w\in\O$ and
$z\in b\O$}\tag2.7$$
(see \cite{2, page~135}).

The kernel $\Lambda(z,w)$ can also be expressed in terms of kernel functions
associated to the boundary (see \cite{4}).
$$\Lambda(w,z)=4\pi L(w,z)^2
+\sum_{i,j=1}^{n-1}\lambda_{ij}L(w,a_i)S(w,a)S(z,a_j)L(z,a),\tag2.8$$
holds for $z,w\in\O$, $z\ne w$.  The coefficients $\lambda_{ij}$
are the same as those appearing in (2.6).  We may express the functions
$S(z,a_j)L(z,a)$ in terms of the other basis $\{{\Cal L}_j\}_{j=1}^{n-1}$
for ${\Cal F}'$ in order to be able to rewrite formula (2.8) in the form
$$\Lambda(w,z)=4\pi L(w,z)^2
+\sum_{i,j=1}^{n-1}\mu_{ij}{\Cal L}_i(z)\,{\Cal L}_j(w)\tag2.9$$
where
${\Cal L}_i(z)=L(z,a_i)S(z,a)$ and, because the matrix
$[\lambda_{ij}]$ is non-singular, so is $[\mu_{ij}]$.

We now suppose that $\O$ is merely an $n$-connected domain in the plane such
that no boundary component of $\O$ is a point.  It is well known that there is
a biholomorphic mapping $\Phi$ mapping $\O$ one-to-one onto a bounded domain
$\O^a$ in the plane with real analytic boundary.  The standard construction
yields a domain $\O^a$ that is a bounded $n$-connected domain with
$C^\infty$ smooth boundary
whose boundary consists of $n$ non-intersecting simple closed real analytic
curves.  Let superscript $a$'s indicate that a kernel function is associated
to $\O^a$.  Kernels without superscripts are associated to $\O$.  The
transformation formula for the Bergman kernels under biholomorphic mappings
gives
$$K(z,w)=\Phi'(z)K^a(\Phi(z),\Phi(w))\overline{\Phi'(w)}.\tag2.10$$
Similarly,
$$\Lambda(z,w)=\Phi'(z)\Lambda^a(\Phi(z),\Phi(w))\Phi'(w).\tag2.11$$
It is well known that the function $\Phi'$ has a single valued
holomorphic square root on $\O$ (see \cite{2, page 43}).
To avoid a discussion of the meaning of the Cauchy transform and
the Szeg\H o projection in non-smooth domains, we shall opt to {\it define\/}
the Szeg\H o and Garabedian kernels associated to $\O$ via
the natural transformation formulas,
$$S(z,w)=\sqrt{\Phi'(z)}\ S^a(\Phi(z),\Phi(w))\overline{\sqrt{\Phi'(w)}},
\tag2.12$$
and
$$L(z,w)=\sqrt{\Phi'(z)}\ L^a(\Phi(z),\Phi(w))\sqrt{\Phi'(w)}.
\tag2.13$$
The Green's functions satisfy
$$G(z,w)=G^a(\Phi(z),\Phi(w))\tag2.14$$
and the functions associated to harmonic measure satisfy
$$\omega_j(z)=\omega_j^a(\Phi(z))\quad\text{ and }\quad
F_j'(z)=\Phi'(z) {F_j^a}'(\Phi(z)).$$
Finally, the Ahlfors map associated to a point $b\in\O$ is defined to be
the solution to the extremal problem, $f_b:\O\to D_1(0)$ with $f_b'(b)>0$ and
maximal.  It is easy to see that the Ahlfors map satisfies
$$f_b(z)=\lambda f_{\Phi(b)}^a(\Phi(z))$$
for some unimodular constant $\lambda$ and it follows that
$f_b(z)$ is a proper holomorphic mapping of $\O$ onto $D_1(0)$.
Furthermore, the transformation formula (2.12) yields that $f_b(z)$
is given by $S(z,b)/L(z,b)$.

It is a routine matter to check that the transformation formulas for the
functions above respect all the formulas given in this section where
the variables range {\it inside\/} the domain.  Therefore, statements
made above about the zeroes of the Szeg\H o kernel function and the
Ahlfors mappings, etc., remain true in this more general setting.

\subhead 3. The Bergman kernel finitely generates itself
\endsubhead
Suppose that $\O$ is an $n$-connected domain in the plane such
that no boundary component of $\O$ is a point.
Because there is a biholomorphic mapping of $\O$ onto a
bounded domain with real analytic boundary, there are plenty of $L^2$
holomorphic functions on $\O$ and so the Bergman projection
can be defined in the standard way as the orthogonal projection
of $L^2(\O)$ onto its subspace of holomorphic functions and the
Bergman kernel is the kernel function associated to this projection.

It is a direct consequence of Theorem~2 in \cite{6} (see also \cite{7})
that any proper holomorphic mapping of $\O$ onto the unit disc
is in a set like $\Cal R_z(\Cal A,N)$ mentioned in \S2.  Because I shall
need some corollaries of the proof given in \cite{6}, I shall give a quick
statement and proof of it here in the simpler case that the proper map
is an Ahlfors map with simple zeroes.

\proclaim{Lemma 3.1}
Suppose $\O$ is an $n$-connected domain domain in the plane such that
no boundary component is a point.  Suppose that $f_a$ is an Ahlfors map
associated to $\O$ such that $f_a$ has exactly $n$ distinct simple
zeroes.  Let $\Cal Z(a)=f_a^{-1}(0)$.  For a positive integer $N$, let
$\Cal R_z(\Cal Z(a),N)$ denote the set of rational combinations of
functions from the set
$$\Cal K(\Cal Z(a),N):=
\{K_m(z,\alpha)\,:\,\alpha\in\Cal Z(a)\text{ and }0\le m\le N\}$$
where $K_0(z,w)$ is the Bergman kernel associated to $\O$ and
$K_m(z,w)$ denotes the function $(\dee^m/\dee\bar w^m)K(z,w)$. The
function $f_a'$ is a linear combination of functions from $\Cal K(\Cal Z(a),0)$
and $f_a'f_a$ is a linear combination of functions from $\Cal K(\Cal Z(a),1)$.
Hence $f_a$ is in $\Cal R_z(\Cal Z(a),1)$.  Furthermore, partial derivatives of
$f_a(z)$ or $f_a'(z)$ with respect to $a$ or $\bar a$ belong to
$\Cal R_z(\Cal Z(a),2)$ as functions of $z$.
\endproclaim

The following theorem is the more general result than Lemma~3.1 that
is proved in \cite{6}.  We state it here for future use.

\proclaim{Theorem 3.2}
Suppose that a domain $\O$ is biholomorphic to a bounded domain in the
plane and that $f$ is a proper holomorphic mapping of $\O$ onto the unit disc.
Let $\Cal Z=f^{-1}(0)$.  There exists a positive integer $N$ such that
$f$ belongs to $\Cal R_z(\Cal Z,N)$.
\endproclaim

A holomorphic function $A(z,w)$ of two complex variables on
an open set in $\C\times\C$ is called {\it algebraic\/}
if there is a holomorphic polynomial $P(a,z,w)$ of three complex variables
such that $A$ satisfies $P(A(z,w),z,w)=0$.  It is a well
established fact that a function $H(z,w)$
which is holomorphic in $z$ and $w$ on a product domain $\O_1\times\O_2$ is
algebraic if and only if, for each fixed $b\in\O_2$, the function $H(z,b)$ is
algebraic in $z$, and for each fixed $a\in\O_1$, the function $H(a,w)$ is
algebraic in $w$ (see Bochner and Martin \cite{9, pages~199-202}).  We shall
say that the Bergman kernel function $K(z,w)$ associated to a domain $\O$ is
{\it algebraic\/} if it can be written as
$R(z,\bar w)$ where $R$ is a holomorphic algebraic function of two
variables on $\{(z,\bar w)\,:\,(z,w)\in\O\times\O\}$.  Because the Bergman
kernel is hermitian, the fact above about separate algebraicity implies that
$K(z,w)$ is algebraic if and
only if, for each point $b\in\O$, the function $K(z,b)$ is an algebraic
function of $z$.  In fact, $K(z,w)$ is algebraic if and only if there
exists a small disc $D_\epsilon(w_0)\subset\O$ such that $K(z,b)$ is
an algebraic function of $z$ for each $b\in D_\epsilon(w_0)$.  We
mention here for future use that if $K(z,w)$ is algebraic, then
differentiating a polynomial identity of the form $P(K(z,w),z,\bar w)=0$
reveals that all the derivatives $(\dee^m/\dee\bar w^m)K(z,w)$ are also
algebraic.

One important consequence of Theorem 3.2 is the following result.

\proclaim{Corollary 3.3}
If $\Phi$ is a biholomorphic mapping between two finitely connected
domains in the plane such that no boundary component is a point
that each have algebraic Bergman kernel functions,
then $\Phi$ must be an algebraic function.
\endproclaim

Corollary 3.3 follows from Theorem~3.2 because if $\Phi:\O_1\to\O_2$
is biholomorphic and if $f_a:\O_2\to D_1(0)$ is an Ahlfors map, then
$f_a\circ\Phi$ is a proper holomorphic map of $\O_1$ onto the unit
disc.  Theorem~3.2 yields that $f_a$ and $f_a\circ\Phi$ are algebraic,
and hence there are germs such that  $\Phi=f_a^{-1}\circ(f_a\circ\Phi)$
and it follows that $\Phi$ is algebraic.

\demo{Proof of Lemma 3.1} Suppose $\O$ and $f_a$ are as in the statement
of the lemma.  It is known that $f_a$ is an $n$-to-one branched covering map
of $\O$ onto $D_1(0)$ (see \cite{2, page~62--70}).  The branch locus
$\Cal B=\{z\in\O: f_a'(z)=0\}$ is a finite set, and for each point
$w_0$ in $D_1(0)-f_a(\Cal B)$, there are exactly $n$ distinct points in
$f_a^{-1}(w_0)$.  Near such a point $w_0$, there is an $\epsilon>0$
such that it is possible to define $n$ distinct holomorphic maps
$\Phi_1(w),\dots,\Phi_n(w)$
on $D_\epsilon(w_0)$ which map into $\O-\Cal B$ such that $f_a(\Phi_k(w))=w$.
These local inverses appear in the following transformation formula for the
Bergman kernels under a proper holomorphic mapping.  Let
$K_D(z,w)=\pi^{-1}(1-z\bar w)^{-2}$ denote the Bergman kernel of the unit
disc (and recall that $K(z,w)$ denotes the Bergman kernel for $\O$).  It
is proved in \cite{2, page~68} that the kernels transform via
$$f_a'(z)K_D(f_a(z),w)=\sum_{k=1}^n K(z,\Phi_k(w))\overline{\Phi_k'(w)}.$$
Since the origin is in $D_1(0)-f(\Cal B)$, we may set $w=0$ in the
transformation formula for the Bergman kernels to obtain
$$f_a'(z)=\pi\sum_{k=1}^n K(z,\Phi_k(0))\overline{\Phi_k'(0)}.\tag3.1$$
This shows that $f_a'(z)$ is a linear combination of functions in
$\Cal K(\Cal Z(a),0)$.  Now differentiate the
transformation formula with respect to $\bar w$ and then set $w=0$
again to obtain
$$2f_a'(z)f_a(z)=\pi\sum_{k=1}^n
\frac{\dee}{\dee\bar w}K(z,\Phi_k(0))\overline{\Phi_k'(0)^2}
+\pi\sum_{k=1}^n K(z,\Phi_k(0))\overline{\Phi_k''(0)}.\tag3.2$$
This shows that $f_a'(z)f_a(z)$ is a linear combination of functions in
$\Cal K(\Cal Z(a),1)$, and so it follows that $f_a(z)$ is in
$\Cal R_z(\Cal Z(a),1)$.

To finish the proof of the lemma, we must analyze the way the
functions $\Phi_k$ depend on $a$ and we shall begin here to write
$\Phi_k(w,a)$ instead of $\Phi_k(w)$.  Let $a_1,a_2,\dots,a_{n-1}$ denote
the simple zeroes of $S(z,a)$ as in \S2 and renumber the $\Phi_k$ so
that $\Phi_k(0,a)=a_k$ for $k=1,\dots,n-1$ and $\Phi_n(0,a)=a$.
If we define $a_n$ to be equal to $a$, then we may write
$\Phi_k(0,a)=a_k$ for all $k$.  Recall that
$f_a(z)=S(z,a)/L(z,a)$.  Hence derivatives of $f_a(z)$, $f_a'(z)$, or
$f_a''(z)$ with
respect to $a$ or $\bar a$ are well defined and easy to compute.
(We shall use the convention that primes denote differentiation with
respect to the $z$ variable.)
Since $\Phi_k(f_a(z),a)=z$ near $z=a_k$, we may differentiate this
formula with respect to $z$ and set $z=a_k$ to see that
$$\Phi_k'(0,a)=1/f_a'(a_k).\tag3.3$$
Next, differentiate the formula
$\Phi_k(f_a(z),a)=z$ twice with respect to $z$ and set $z=a_k$ to obtain
$\Phi_k''(0,a)f_a'(a_k)^2+\Phi_k'(0,a)f_a''(a_k)=0$.
This shows that
$$\Phi_k''(0,a)=-f_a''(a_k)/f_a'(a_k)^3.\tag3.4$$
To finish the proof, we must consider the way in which the zeroes $a_k$
depend on $a$.  We now write $a_k(a)$ in order to regard $a_k$ as a function
of $a$.  (Of course, $a_n(a)=a$ and it is only for $k=1,\dots,n-1$ that
we need to exert ourselves.)  Let $A_0$ be
a fixed point in $\O$ such that the zeroes of $S(z,A_0)$ are simple.
Since the points $a_k(A_0)$ are distinct, we may choose an
$\epsilon>0$ such that the closed discs of radius $\epsilon$ about
the points $a_k(A_0)$  are mutually disjoint and each contained in
$\O$.  Thus, $a_k(A_0)$ is the only zero of $S(z,A_0)$ in the closure
of $D_\epsilon(a_k(A_0))$.  The dependence
of the zeroes of $S(z,a)$ on $a$ can be read off from the formula,
$$a_k(a)=\frac{1}{2\pi i}\int_{|z-a_k(A_0)|=\epsilon}
z\ \frac{\frac{\dee}{\dee z}S(z,a)}{S(z,a)}\ dz,$$
which is valid when $a$ is close to $A_0$.
Because $S(z,a)$ is antiholomorphic in $a$, this formula shows that
$a_k(a)$ {\it is an antiholomorphic function of\/} $a$ near $A_0$.

We may now differentiate (3.1) and (3.2) with respect to $a$ or $\bar a$
and use the complex chain rule together with (3.3) and (3.4) to complete
the proof of the lemma.
\enddemo

We now turn to the proof of Theorem~2.1.
Assume that $\O$ is an $n$-connected domain in the plane such
that no boundary component of $\O$.  Lemma~3.1 yields that the Ahlfors maps
associated to $\O$ are rational combinations of the type
$\Cal R_z(\Cal Z(a),1)$.

The zeroes of $S(z,a)$ in the $z$ variable become simple zeroes
as $a$ approaches the boundary (\cite{3}) and so for $a\in\O$
that are sufficiently close to the boundary we know that $S(z,a)$
has exactly $n-1$ simple zeroes in the $z$ variable.  As in \S2, we shall
denote these zeroes by $a_i$, $i=1,\dots,n-1$ and we shall let
$\Cal Z(a)$ denote the set of $n$ points, $\{a,a_1,\dots,a_{n-1}\}$.
Because the Ahlfors map $f_a$ is a proper holomorphic map onto
the unit disc, and $a$ together with $a_1,\dots,a_{n-1}$ are the
zeroes of $f_a$, it follows that the classical Green's function $G(z,w)$
associated to $\O$ satisfies
$$\frac{1}{2}\ln |f_a(z)|^2=G(z,a)+\sum_{i=1}^{n-1}G(z,a_i).\tag3.5$$

We now differentiate (3.5) with respect to $z$ to obtain
$$\frac{f_a'(z)}{2f_a(z)}=\frac{\dee}{\dee z} G(z,a)+
\sum_{i=1}^{n-1} \frac{\dee}{\dee z} G(z,a_i).$$
During the course of the proof of Lemma~3.1, we showed that
the zeroes $a_1,\dots,a_{n-1}$ are antiholomorphic functions of $a$
when $a$ is near the boundary.
Next, we differentiate with respect to $a$ and use the complex chain
rule to obtain
$$\frac{\dee}{\dee a}\left(\frac{f_a'(z)}{2f_a(z)}\right)
=\frac{\dee^2 G(z,a)}{\dee z\dee a}+
\sum_{i=1}^{n-1} \frac{\dee^2 G(z,a_i)}{\dee z\dee\bar a_i}
\frac{\dee\bar a_i}{\dee a}.\tag3.6$$
Let $R(z,a)$ denote the left hand side of (3.6).
Lemma~3.1 yields that, as a function of $z$, $R(z,a)$ is in
$\Cal R_z(\Cal Z(a),2)$.
The function on the right hand side of (3.6) can be rewritten to yield
$$R(z,a)=-\frac{\pi}{2}\Lambda(z,a)-\frac{\pi}{2}
\sum_{i=1}^{n-1}K(z,a_i)\frac{\dee\bar a_i}{\dee a}.$$
This last formula shows that, for each fixed $a$ sufficiently close
to the boundary of $\O$, the function $\Lambda(z,a)$ is in
$\Cal R_z(\Cal Z(a),2)$.  Note that we may also state that
$\Lambda(z,a_k)$ is in $\Cal R_z(\Cal Z(a_k),2)$ for $k=1,\dots,n-1$
when $a$ is sufficiently close to the boundary.

Solve formula (2.6) for $S(z,w)^2$ and formula (2.9)
for $L(z,w)^2$ and divide the two, noting that the Ahlfors map
$f_w(z)$ is equal to $S(z,w)/L(z,w)$,
to obtain
$$f_w(z)^2=\frac{K(z,w)-\sum_{i,j=1}^{n-1}
\lambda_{ij}{\Cal L}_i(z)\,\overline{{\Cal L}_j(w)}}
{\Lambda(w,z)-
\sum_{i,j=1}^{n-1}\mu_{ij}{\Cal L}_i(z)\,{\Cal L}_j(w)}.\tag3.7$$
Recall that $\Cal L_i(z):=L(z,a_i)S(z,a)$.  We shall use formula (3.7)
to show that the functions $\{{\Cal L}_j(z)\}_{j=1}^{n-1}$ are in
$\Cal R_z(\Cal A(a),2)$ where
$$\Cal A(a)=\Cal Z(a)\cup_{k=1}^{n-1}\Cal Z(a_k)$$
is the set described before the statement of Theorem~2.1.
We may manipulate (3.7) to obtain
$$\sum_{i,j=1}^{n-1}{\Cal L}_i(z)
\left(\lambda_{ij}\,\overline{{\Cal L}_j(w)}-
\mu_{ij}f_w(z)^2{\Cal L}_j(w)\right)=
K(z,w)-f_w(z)^2\Lambda(w,z).\tag3.8$$
Our plan now is to plug into this formula $w=a_1,a_2,\dots,a_{n-1}$ where
the $a_k$'s are the zeroes of the Szeg\H o kernel $S(z,a)$ that appear in the
definition of the functions ${\Cal L}_i(z)$.  We shall obtain a system of
$n-1$ equations that we can use to solve for the
functions ${\Cal L}_i(z)$ and thereby see that the
${\Cal L}_i(z)$ are in $\Cal R_z(\Cal A(a),2)$.
Define
$$A_{ik}(z)=
\sum_{j=1}^{n-1}
\left(\lambda_{ij}\,\overline{{\Cal L}_j(a_k)}-
\mu_{ij}f_{a_k}(z)^2{\Cal L}_j(a_k)\right).$$
Notice that, because ${\Cal L}_j(z)=L(z,a_j)S(z,a)$ and because $S(a_k,a)=0$,
we see that 
$${\Cal L}_j(a_k)=
\cases 0 &\text{if }j\ne k \\
q_k &\text{if }j=k
\endcases
$$
where, because the zeroes of $S(z,a)$ are {\it simple\/} zeroes, $q_k$ is
a non-zero number given by $V'(a_k)$ where $V(z)=S(z,a)/(2\pi)$.
This shows that the $(n-1)\times(n-1)$ matrix $[{\Cal L}_j(a_k)]$ is
non-singular.  Notice that
$$A_{ik}(a)=
\sum_{j=1}^{n-1}\lambda_{ij}\,\overline{{\Cal L}_j(a_k)}$$
because $f_{a_k}(a)$ is zero via (2.3) and the fact that
$S(a,a_k)=0$.  It therefore follows from the matrix identity
$$[A_{ik}(a)]=[\lambda_{ij}][\ \overline{{\Cal L}_j(a_k)}\ ]$$
and Hejhal's theorem
that $[A_{ik}(a)]$ is non-singular.  Hence, for $z$ in a neighborhood
of $a$, the matrix $[A_{ik}(z)]$ will be non-singular and we can
use Cramer's rule to solve the $n-1$ equations obtained from (3.8) by
plugging in $w=a_1,a_2,\dots,a_{n-1}$ to see that each of the
functions ${\Cal L}_j(z)$ is a rational combination of the functions
$f_{a_k}(z)$ and $K(z,a_k)-f_{a_k}(z)^2\Lambda(a_k,z)$,
$k=1,\dots,n-1$.  Since all of these functions are in $\Cal R_z(\Cal A(a),2)$,
we conclude that ${\Cal L}_j(z)$ is in
$\Cal R_z(\Cal A(a),2)$.

It now follows from (2.9) that $L(z,a)^2$ is also in
$\Cal R_z(\Cal A(a),2)$.  Next, we multiply (2.4) by
$L(z,a)\overline{L(w,a)}$ to see that
$$
S(z,w)L(z,a)\overline{L(w,a)}$$
is a rational combination of the functions $f_a(z)$ and $\overline{f_a(w)}$
times a linear combination of the functions
$S(z,a)L(z,a)$, $\overline{S(w,a)L(w,a)}$, and
$S(z,a_i)L(z,a)$, $\overline{S(w,a_j)L(w,a)}$  for $i,j=1,\dots,n-1$.
This shows that $S(z,w)L(z,a)\overline{L(w,a)}$ is in
$\Cal R_{z,\bar w}(\Cal A(a),2)$ because $f_a$ is in $\Cal R_z(\Cal A(a),2)$,
and $S(z,a)L(z,a)=f_a(z)L(z,a)^2$ is in $\Cal R_z(\Cal A(a),2)$ since
$f_a$ and $L(z,a)^2$ are, and each of the functions
$S(z,a_i)L(z,a)$ is in $\Cal R_z(\Cal A(a),2)$ because, as mentioned
in \S2, the linear span of $\{\Cal L_i:i=1,\dots,n-1\}$ is the same
as the linear span of $\{S(z,a_i)L(z,a):i=1,\dots,n-1\}$.

Finally, we multiply (2.6) by $L(z,a)^2\overline{L(w,a)^2}$
to obtain
$$\gather K(z,w)L(z,a)^2\overline{L(w,a)^2}= \\
4\pi [S(z,w) L(z,a)\overline{L(w,a)}]^2+
L(z,a)^2\overline{L(w,a)^2}\sum_{i,j=1}^{n-1}
\lambda_{ij}{\Cal L}_i(z)\,\overline{{\Cal L}_j(w)}.
\endgather$$
We have shown that the right hand side of this equation is composed of
functions in $\Cal R_{z,\bar w}(\Cal A(a),2)$.  Since we know that
$L(z,a)^2$ is in $\Cal R_z(\Cal A(a),2)$, we may divide the equation
by $L(z,a)^2\overline{L(w,a)^2}$ to see that $K(z,w)$ is in
$\Cal R_{z,\bar w}(\Cal A(a),2)$.  Now formula (2.6) shows that
$S(z,w)^2$ is in $\Cal R_{z,\bar w}(\Cal A(a),2)$.
Similar reasoning using (2.5) and (2.9) reveals that $\Lambda(z,w)$
and $L(z,w)^2$ are in $\Cal R_{z,w}(\Cal A(a),2)$.
The proof of Theorem~2.1 is now complete.

\subhead 4. The double of a domain, proper maps, and the Bergman kernel
\endsubhead
It is well known that the Bergman kernel can be extended to the double
of a smooth finitely connected domain as a meromorphic {\it
differential}.  In this section, we show that certain simple
combinations of the Bergman kernel and a proper holomorphic map onto
the unit disc extend to the double as meromorphic {\it functions}.

\proclaim{Theorem 4.1}
Suppose that $\O$ is a bounded $n$-connected domain in the plane
bounded by $n$ non-intersecting $C^\infty$ smooth real analytic curves and
suppose that $f:\O\to D_1(0)$ is a proper holomorphic map.  Then
$f(z)$ extends meromorphically to the double of $\O$.  Also,
$K(z,w)/f'(z)$ extends meromorphically to the double of $\O$ as a
function of $z$ for each $w\in\Obar$ and so do the functions
$K_m(z,w)/f'(z)$ for each non-negative integer $m$ (where
$K_0(z,w)$ denotes the Bergman kernel of $\O$ and
$K_m(z,w)=(\dee^m/\dee\bar w^m)K(z,w)$).
Furthermore, if the zeroes of $f(z)$ in $\O$ are simple
zeroes, then for all but possibly finitely many $b$ in $\Obar$, the
functions $f(z)$ and $K(z,b)/f'(z)$, as functions of $z$, generate
the field of meromorphic functions on the double of $\O$, i.e.,
they form a primitive pair for the double of $\O$  as in \cite{10, page~249}.
Hence, there is an irreducible polynomial $P(z,w)$ such that
$P(f(z),K(z,b)/f'(z))\equiv0$ on $\O$, and for fixed $b$ as above, we
may state that $K(z,b)/f'(z)$ is an algebraic function of $f(z)$.
\endproclaim

We shall be able to combine Theorems 2.1 and 4.1 to obtain the following
theorem.

\proclaim{Theorem 4.2}
Suppose that $\O$ is a finitely connected domain in the plane such
that no boundary component is a point and suppose that
$f:\O\to D_1(0)$ is a proper holomorphic map.
There is a point $b\in\O$ such that the Bergman kernel $K(z,w)$
associated to $\O$ is a rational combination of the six functions
$f(z)$, $f'(z)$, and $K(z,b)$ and the conjugates of $f(w)$, $f'(w)$, and
$K(w,b)$.  Furthermore, there is an irreducible polynomial $P(z,w)$ such that
$P(f(z),K(z,b)/f'(z))\equiv0$ on $\O$, and therefore $K(z,b)/f'(z)$ is an
algebraic function of $f(z)$.  Hence, $K(z,w)$ is an algebraic function of
$f(z)$ and $f'(z)$ and the conjugates of $f(w)$ and $f'(w)$.
\endproclaim

Theorem 4.2 gives a rather clean answer to a question posed in
\cite{5}, ``When is the Bergman kernel associated to a domain algebraic?''

\proclaim{Corollary 4.3}
Suppose that $\O$ is an $n$-connected domain in the plane such that no
boundary component is a point.  The Bergman kernel associated to $\O$ is
algebraic if and only if there
exists a proper holomorphic mapping of $\O$ onto the unit disc which is
algebraic, i.e., if and only if there is a non-constant algebraic function
$f$ which has a single valued holomorphic sheet over $\O$ such that $|f(z)|$
tends to one as $z$ tends to the boundary of $\O$.  Furthermore, if the
Bergman kernel is algebraic, then every proper holomorphic mapping of
$\O$ onto the unit disc is algebraic, the functions $F_1',\dots,F_{n-1}$
are algebraic, and so are the Szeg\H o kernel, the Garabedian kernel, and
the kernel $\Lambda(z,w)$.
\endproclaim

This corollary implies, for example, that all the classical domain
functions associated to the $2$-connected domain given by
$\{z\,:\, |z+1/z|<r\}$ are algebraic when $r$ is a real constant
bigger than $2$.  An interesting application of our results is
that every $2$-connected domain is biholomorphic to a domain with
algebraic Bergman kernel.  Indeed, standard arguments (see \cite{12,
pages 10-11}) show that the modulus of $\{z\,:\, |z+1/z|<r\}$ 
is a continuous increasing function of $r$ which goes to zero as $r$
approaches $2$ from above (since the domains have a narrow ``pinch''
near $\pm 1$ as $r\to2$) and which goes to infinity as $r$ tends to
infinity (since ``fat'' annuli can be put inside the domains).
Hence, each $2$-connected domain is biholomorphic to exactly one domain
of the form $\{z\,:\, |z+1/z|<r\}$ with $r>2$.  Since Corollary~3.3
states that biholomorphic maps between domains with algebraic
Bergman kernel functions are algebraic, we may think of
$\{z\,:\, |z+1/z|<r\}$ as being the defining member of a class.
This domain also has the virtue that the mapping $f(z)=(1/r)(z+1/z)$
is a $2$-to-one branched cover of $A(r)$ onto the unit disc that
extends to be a one-to-one biholomorphic mapping from each connected
component of the complement of $A(r)$ in the Riemann sphere onto
the complement of the unit disc in the Riemann sphere.

An interesting problem remains to extend these ideas to
$n$-connected domains.  Can every $n$-connected domain be mapped
to a domain of the form
$$\{ |z+\sum_{k=1}^{n-1}a_k/(z-b_k)|<r\}?$$

When Corollary~4.3 is combined with results of \S3,
we obtain the following result.

\proclaim{Theorem 4.4}
Suppose $\O$ is a finitely connected domain in the plane such that
no boundary component is a point.  The following conditions are
equivalent.
{\roster
\item The Bergman kernel associated to $\O$ is algebraic.
\item The Szeg\H o kernel associated to $\O$ is algebraic.
\item There exists a single proper holomorphic mapping of $\O$ onto
the unit disc which is algebraic.
\item Every proper holomorphic mapping of $\O$ onto
the unit disc is algebraic.
\endroster}
\endproclaim

We remark that results in \cite{4} giving explicit formulas for
the Poisson kernel in terms of the Szeg\H o kernel also allow us to
deduce that if the Szeg\H o kernel associated to a finitely connected
domain in the plane such that no boundary component is a point is
algebraic, then the Poisson kernel of the domain is a real algebraic
function of the variable that runs over the boundary.  (It is possible
to formulate a converse of this statement because the Poisson kernel
is related to a derivative $(\dee/\dee\bar w)G(z,w)$ of the Green's
function and the Bergman kernel is a constant times
$(\dee^2/\dee z\dee\bar w)G(z,w)$.  However, because the statement is
rather awkward, and because the interesting direction of the
implication is the one we stated here, we omit it.)

Whenever a proper map $f$ of a domain onto the unit disc is algebraic, then
so are the Schwarz Reflection Functions associated to the boundary curves
of the domain because the antiholomorphic reflection functions can be
written as $R(z)=f^{-1}(1/\,\overline{f(z)}\,)$ near the boundary.
Hence, the techniques used in this paper might have applications
to the sorts of questions studied in Shapiro \cite{13}.

\demo{Proof of Theorem 4.1}
Suppose that $\O$ and $f$ are as in the statement of Theorem~4.1.  It is
well known that $f$ extends holomorphically past the boundary of $\O$,
that $f'$ is non-vanishing on $b\O$, and that there is a positive
integer $m$ such that $f$ is an $m$-to-one branched covering map (see
\cite{2, page~62-66}).  Furthermore, $|f(z)|=1$ for $z\in b\O$.

Let $\Oh$ denote the double of $\O$ and let $R(z)$ denote
the antiholomorphic involution on $\Oh$ which fixes the boundary of
$\O$.  Let $\Ot=R(\O)$ denote the reflection of $\O$ across the
boundary.  It is easy to see that $f$ extends to be a meromorphic function
on $\Oh$ because $f(z)=1/\,\overline{f(z)}$ for $z\in b\O$ and, since
$R(z)=z$ on $b\O$, it follows that
$$f(z)=1/\,\overline{f(R(z))}\text{\quad for }z\in b\O.$$
The function on the left hand side of this formula is holomorphic
on $\O$ and the function on the right hand side is meromorphic
on $\Ot$ and the two functions extend continuously to $b\O$ from
opposite sides and agree on $b\O$.  Hence, the function given
by $f(z)$ on $\Obar$ and $1/\,\overline{f(R(z))}$ on $\Ot$ is
meromorphic on $\Oh$.  We next note that, because $\ln|f(z)|^2=0$
for $z\in b\O$, we may differentiate $\ln|f(z(t))|^2$ with respect
to $t$ when $z(t)$ parameterizes the boundary to obtain
$$\frac{f'(z(t))}{f(z(t))}z'(t)
+\overline{f'(z(t))z'(t)/f(z(t))}=0.$$
Dividing this equation by $|z'(t)|$ reveals that
$$\frac{f'(z)}{f(z)}T(z)=-\overline{f'(z)T(z)/f(z)},\qquad
\text{for }z\in b\O.\tag4.1$$
Notice the similarity of this formula to (2.7).  Let $w$ be a fixed
point in $\Obar$.  The conjugate of (2.7) is
$$K(z,w)T(z)=-\overline{\Lambda(w,z)}\,\overline{T(z)}\qquad\text{for
$z\in b\O$}.\tag4.2$$
We may divide (4.2) by (4.1) to obtain
$$\frac{f(z)K(z,w)}{f'(z)}=\overline{f(z)\Lambda(w,z)/f'(z)}\qquad\text{for
$z\in b\O$}.$$
Finally, we may replace $z$ by $R(z)$ (because $R(z)=z$ on $b\O$) in the
right hand side of this last equation to obtain
$$\frac{f(z)K(z,w)}{f'(z)}=\overline{f(R(z))\Lambda(w,R(z))/f'(R(z))}
\qquad\text{for $z\in b\O$}.\tag4.3$$
The function on the left hand side of (4.3) is meromorphic on $\O$ and
extends holomorphically up to the boundary.  The function on the right
hand side is meromorphic on $\Ot$ and extends to the boundary of $\O$
from the outside of $\O$.  Since these two agree on the boundary,
the right hand side defines the meromorphic extension of $f(z)K(z,w)/f'(z)$
to $\Oh$.  Since $f(z)$ also extends to be meromorphic on $\Oh$, we may
divide by $f(z)$ to see that $K(z,w)/f'(z)$ extends meromorphically
to $\Oh$ for each $w\in\Obar$.

This last argument can be repeated with $K_m(z,w)$ and
$(\dee^m/\dee w^m)\Lambda(z,w):= \Lambda_m(z,w)$ in place of $K(z,w)$
and $\Lambda(z,w)$, respectively, because (2.7) can be differentiated
with respect to $w$ and then conjugated to yield
$$K_m(z,w)T(z)=-\overline{\Lambda_m(w,z)}\,\overline{T(z)}\qquad\text{for
$z\in b\O$}.$$
Hence, $K_m(z,w)/f'(z)$ extends meromorphically to $\Oh$ for each $w$
in $\Obar$.

To finish the proof of Theorem~4.1, we must study conditions under which
$f(z)$ and $K(z,w)/f'(z)$ form a {\it primitive pair}.  We now assume
that the zeroes of $f(z)$ are simple zeroes.  Let $a_1,\dots,a_m$
denote the $m$ distinct zeroes of $f$ on $\O$.  It is easy to see that
the order of $f$ as a meromorphic function on $\Oh$ is $m$ (because the
$m$ zeroes of $f$ on $\O$ get reflected by $R$ to $m$ poles on $\Ot$).
Hence, to prove that $f(z)$ and $K(z,w)/f'(z)$ form a primitive pair for most
$w$, we must show that $K(z,w)/f'(z)$ separates the $m$ points
$a_1,\dots,a_m$ as a function of $z$ for most $w$ (see \cite{1, page~321-324}).
Let $S_{ij}$ denote the set of points $w$ in $\Obar$ such that
$K(a_i,w)=c_{ij}K(a_j,w)$ where $c_{ij}$ is a non-zero constant given by
$f'(a_i)/f'(a_j)$.
It is easy to see that the linear span of the set of functions of $z$
given by $\Cal K:=\{K(z,w)\,:\,w\in\O\}$ is dense in $H^2(\O)$.  (Indeed,
if $h(z)$ is a function in $H^2(\O)$ that is orthogonal to all the
functions in the spanning set, then $h(w)=0$ for all $w\in\O$ by the
reproducing property of the Bergman kernel, and this proves the
density.)  Let $P$ be a polynomial of one complex variable such
that $P(a_i)\ne c_{ij}P(a_j)$.  Since it is possible to choose a
linear combination of functions in the spanning set $\Cal K$ which
converge uniformly on compact subsets of $\O$ to $P$, it follows
that the function of $w$ given by $K(a_i,w)-c_{ij}K(a_j,w)$ cannot be
identically zero on $\O$.  Since this function extends to be
holomorphic on a neighborhood of $\Obar$ as a function of $w$, it
follows that $S_{ij}$ is a finite set.  Hence, the set of points
$\cup_{i<j}S_{ij}$ where $K(z,w)/f'(z)$ might fail to separate the
$m$ points $a_1,\dots,a_m$ as a function of $z$ is at most a finite
subset of $\Obar$.

Finally, it is a classical fact that, whenever $g$ and $h$ form a primitive
pair on a compact Riemann surface, there exists an irreducible
polynomial $P(z,w)$ such that $P(g(z),h(z))\equiv0$ on the surface
(see Farkas and Kra \cite{10, page~248-250}).
This finishes the proof of Theorem~4.1.
\enddemo

\demo{Proof of Theorem 4.2}
Suppose that $\O$ is a domain as described in the hypotheses of
Theorem~4.2.  Suppose for the moment that the zeroes of $f$ on $\O$
are simple zeroes.  Theorem~2.1 yields that the Bergman kernel
$K(z,w)$ generates itself.
Let  $\Phi$ denote a biholomorphic mapping which maps
$\O$ one-to-one onto a bounded domain $\O^a$ in the plane with real
analytic boundary curves and let $\phi=\Phi^{-1}$ and $F=f\circ \Phi^{-1}$.
Let $K^a(z,w)$ denote the Bergman kernel
associated to $\O^a$.  Since $F$ is a proper holomorphic mapping of
$\O^a$ onto the unit disc with simple zeroes, Theorem~4.1 yields that
there exists a point
$b\in\O^a$ such that $F(z)$ and $K^a(z,b)/F'(z)$ form a primitive
pair for the field of meromorphic functions on the double of $\O^a$.
Let $B=\phi(b)$.

The transformation formula (2.10) for the Bergman kernels
under biholomorphic mappings can also be written in the form
$$\phi'(z)K(\phi(z),w)=K^a(z,\Phi(w))\overline{\Phi'(w)}.$$
Hence, we may transform an expression of the form
$$K(z,w)/f'(z)$$
by replacing $z$ by $\phi(z)$ and by multiplying the whole thing by
$\phi'(z)/\phi'(z)$ to obtain
$$K(\phi(z),w)/f'(\phi(z))= \overline{\Phi'(w)}K^a(z,\Phi(w))/F'(z).\tag4.4$$
Similar reasoning shows that
$$K(z,w)/f'(z)= \overline{\Phi'(w)}K^a(\Phi(z),\Phi(w))/F'(\Phi(z)).\tag4.5$$
Differentiating (4.4) repeatedly with respect to $\bar w$
yields that the linear span of the set
$$\{K_m(\phi(z),w)/f'(\phi(z)): m\le N\}$$
is equal to the linear span of
$$\{K_m^a(z,\Phi(w))/F'(z): m\le N\}.$$
Theorem~4.1 states that functions of the form $K_m^a(z,\Phi(w))/F'(z)$
are rational combinations of $F(z)$ and $K^a(z,b)/F'(z)$.
It follows that the field of functions generated by
functions of the form $K_m(\phi(z),w)/f'(\phi(z))$ is equal to the field of
functions generated by $F(z)$ and $K^a(z,b)/F'(z)$.  Finally, by
replacing $z$ by $\Phi(z)$ and by using (4.5), we see that
the field of functions generated by
functions of the form $K_m(z,w)/f'(z)$ is equal to the field generated
by $f(z)=F(\Phi(z))$ and $K(z,B)/f'(z)$.  It now follows from Theorem~2.1 that
the Bergman kernel $K(z,w)$ is in the field generated by
by $f(z)$, $f'(z)$, and $K(z,B)$, and conjugates of
$f(w)$, $f'(w)$, and $K(w,B)$.

Similar reasoning using (4.5) shows that, by replacing $z$ by $\Phi(z)$,
a polynomial identity of the form $P^a(F(z),K^a(z,b)/F'(z))=0$ can be
transformed into one of the form $P(f(z),K(z,B)/f'(z))=0$ and the proof
is complete in case $f$ has simple zeroes.

To finish the proof, we must treat the case of a proper holomorphic
mapping $f:\O\to D_1(0)$ that does not have simple zeroes.  This is
easy, however, because we may compose $f$ with a M\"obius
transformation $\psi$ so that $\psi\circ f$ has simple zeroes.  Since
$\psi$ is rational, it is easy to see that $\psi\circ f$ and
$(\psi\circ f)'$ are contained in the function field generated by
$f$ and $f'$ and the first part of the theorem follows.  An explicit
calculation shows that a polynomial identity of the form
$P((\psi\circ f)(z),K(z,B)/(\psi\circ f)'(z))=0$
can be converted to one of the form
$\widetilde P(f(z),K(z,B)/f'(z))=0$ after clearing the denominator terms
by multiplying by a suitable polynomial in $f(z)$.
The proof of Theorem~4.2 is complete.
\enddemo

We remark that if we define $I(z,w)=K(z,w)/[f'(z)\overline{f'(w)}]$ and
$I^a(z,w)=K^a(z,w)/[F'(z)\overline{F'(w)}]$ in the proof of Theorem~4.2
above, we obtain a rather interesting invariant.  The transformation
formula (2.10) for the Bergman kernels together with the fact that
$f'(z)=F'(\Phi(z))\Phi'(z)$ yields that
$$I(z,w)=I^a(\Phi(z),\Phi(w)),$$
and polynomials satisfying $P(f(z),I(z,b))=0$ can be viewed as genuine
algebraic geometric invariants.

\subhead 5. The case of algebraic kernel functions
\endsubhead
When the results of this paper are combined with results proved in \cite{5},
we can show that if the Bergman kernel associated to a finitely connected
domain in the plane is algebraic, then the Bergman kernel is a
rational combination of just {\it two\/} holomorphic functions of one
variable and it is an algebraic function of a single proper
holomorphic map.  To prove this, we shall need to use the following
result from \cite{5}.

\proclaim{Theorem 5.1}
Suppose that $\O$ is a finitely connected domain in the plane such that no
boundary component is a point.  If the Bergman or the Szeg\H o kernel
associated to $\O$ is algebraic, then $\O$ can be realized
as a subdomain of a compact Riemann surface $\widetilde{\Cal R}$
such that all the kernel functions
$S(z,w)$, $L(z,w)$, $K(z,w)$, $\Lambda(z,w)$ extend to
$\widetilde{\Cal R}\times\widetilde{\Cal R}$
as single valued meromorphic functions.  Furthermore, the Ahlfors
maps $f_a(z)$ and every proper holomorphic mapping from $\O$ to
the unit disc extend to be single valued meromorphic functions on
$\widetilde{\Cal R}$.  Also, the functions $F_k'(z)$, $k=1,\dots,n-1$,
extend to be single valued meromorphic functions on 
$\widetilde{\Cal R}$.  Furthermore, the complement of $\O$ in
$\widetilde{\Cal R}$ is connected.
\endproclaim

When this theorem is combined with results in \S4, we derive the
following result.

\proclaim{Theorem 5.2}
Suppose $\O$ is a finitely connected domain in the plane such that no
boundary component is a point and suppose that
the Bergman kernel $K(z,w)$ associated to $\O$ is algebraic.
If $f:\O\to D_1(0)$ is a proper holomorphic mapping, then
$K(z,w)$ is an algebraic function of $f(z)$ and $\overline{f(w)}$ and
there exists an irreducible polynomial $P(u,v,w)$ of three complex
variables such that
$$P(K(z,w),f(z),\,\overline{f(w)}\,)\equiv0.$$
Furthermore, there exists a holomorphic function $g(z)$ on $\O$ such
that $K(z,w)$ is a rational combination of $f(z)$, $g(z)$, and conjugates
of $f(w)$ and $g(w)$.
\endproclaim

Let $\Oh$ denote the double of $\O$ and let $R(z)$ denote the
antiholomorphic reflection function associated to the construction of
$\Oh$.  The polynomial equation in
Theorem~5.2 reveals that it is possible to analytically continue $K(z,w)$
to $\Oh\times\Oh$ as a finitely valued multivalued function
with algebraic singularities.  Indeed, since
$f(z)=\overline{f(R(z))^{-1}}$ for $z\in b\O$, the expression
$P(K(z,w),\overline{f(R(z))^{-1}},\overline{f(w)})$ is equal to
$P(K(z,w),f(z),\,\overline{f(w)}\,)$ for $z\in b\O$ and it
defines the holomorphic continuation of $K(z,w)$ in the $z$ variable to the
reflection $R(\O)$ for each fixed $w\in\O$.  Now the same thing can be done
in the $w$ variable for each fixed $z\in\Oh$ to complete the extension to
$\Oh\times\Oh$.

\demo{Proof of Theorem 5.2}
Theorem~5.1 states that the proper holomorphic map $f$ extends as a
meromorphic function to the compact Riemann surface $\widetilde{\Cal R}$.
Suppose that the order of $f$ on $\widetilde{\Cal R}$ is $m$.  Choose
a point $\lambda\in\C$ with $|\lambda|>1$ so that $f^{-1}(\lambda)$
consists of $m$ distinct points.  We may construct a meromorphic
function $g$ on $\widetilde{\Cal R}$ as in Farkas and Kra \cite{10,
page 248-249} which is holomorphic on $\O\subset\widetilde{\Cal R}$
such that $f$ and $g$ form a primitive pair for the field of
meromorphic functions on $\widetilde{\Cal R}$.  The construction of
$\widetilde{\Cal R}$ in \cite{5} reveals
that functions of the form $K(z,b)$ and $K_m(z,b)$ extend to
$\widetilde{\Cal R}$ as meromorphic functions of $z$.  Hence, by the
properties of a primitive pair, these functions are generated by $f$
and $g$.  Now Theorem~4.2 yields that $K(z,w)$ is a rational combination
of $f(z)$, $g(z)$, and conjugates of $f(w)$ and $g(w)$.  Finally,
there is an irreducible polynomial $G(z,w)$ such that $G(f(z),g(z))\equiv0$
on $\widetilde{\Cal R}$.  It follows that $g$ is an algebraic function
of $f$, and hence that $K(z,w)$ is an algebraic function of
$f(z)$ and the conjugate of $f(w)$.  Hence, a polynomial identity of
the form $P(K(z,w)f(z),\overline{f(w)})\equiv0$ as in the statement of
the theorem holds.  This completes the proof.
\enddemo

\enddocument